\documentstyle[12pt]{article}
\newcommand{\too}{\longrightarrow}
\newcommand{\om}{\omega}
\newcommand{\an}{{\pi_{\#}}}

\newcommand{\X}{{\cal X}}
\newcommand{\G}{{\cal G}}
\newcommand{\F}{{\cal F}}

\newcommand{\D}{{\cal D}}

\newcommand{\p}{{\cal P}}

\newcommand{\di}{\displaystyle}
\newcommand{\Om}{\Omega}
\newcommand{\na}{\nabla}
\newcommand{\wi}{\widetilde}
\newcommand{\al}{\alpha}
\newcommand{\be}{\beta}
\newcommand{\ga}{\gamma}
\newcommand{\Ga}{\Gamma}

\newcommand{\De}{\Delta}
\newcommand{\de}{\delta}

\def \reel{ {\rm I}\!{\rm R} }

 \def \nat{ { {\rm I}\!{\rm N}} }

 \def \rat{ {\rm Q}\kern-.65em {}^{{}_/ }}

\newtheorem{Def}{Definition}[section]
\newtheorem{th}{Theorem}[section]
\newtheorem{pr}{Proposition}[section]

\newtheorem{co}{Corollary}[section]
\title{Killing-Poisson tensors on Riemannian manifolds
 } \author{M. Boucetta\\
Facult\'e des Sciences et Techniques \\
BP 549 Marrakech\\
Morocco
\\
Email: {\it boucetta@fstg-marrakech.ac.ma
}
\footnote{Recherche men\'ee dans le cadre du Programme
 Th\'ematique d'Appui \`a la Recherche Scientifique PROTARS III.}} \date{ }
\begin{document}
\maketitle
\parindent=0cm

{\bf Abstract.} We  introduce a new class of Poisson structures on
a Riemannian manifold. A Poisson structure in this class will be
called a Killing-Poisson structure. The class of Killing-Poisson
structures contains the class of symplectic structures, the class
of Poisson structures studied in{ \it  (Differential Geometry and
its Applications, {\bf Vol. 20, Issue 3} (2004), 279--291)} and
the class of Poisson structures induced by some infinitesimal Lie
algebras actions on Riemannian manifolds.  We show that some
classical results on symplectic manifolds (the integrability of
the Lie algebroid structure associated to a symplectic structure,
the non exactness of a symplectic structure on a compact manifold)
remain valid for regular Killing-Poisson structures.

\bigskip



\section{Introduction and main results}
In the present paper, we pursue our investigations on the
interactions between Riemannian geometry and Poisson geometry
initialized in $[1]$, $[2]$ and $[3]$, by introducing a new class
of Poisson structures. A Poisson structure in this new class will
be called a {\bf Killing-Poisson structure}. We show that this
terminology is appropriate by pointing out that there is a notion
of Killing multi-vector fields as a generalization of Killing
vector fields (cf. Proposition 1.1) and, by definition, a
Killing-Poisson structure is a Poisson structure whose associated
bivector field is a Killing bivector field. We show that the class
of Killing-Poisson structures contains the class of symplectic
structures, the class of Poisson
 structures studied in $[2]$ and the class
 of Poisson
 structures  induced by some infinitesimal Lie algebras actions on Riemannian
  manifolds.
  We show
 that some classical results on symplectic manifolds
  (the integrability of the Lie algebroid
 structure associated to a symplectic structure,
 the non exactness of a symplectic
 structure
 on a compact manifold) remain valid for regular Killing-Poisson structures.

To state our results, we first explain some fundamental notions
which we use in this paper.

Recall that a Poisson structure  on a manifold $M$ is an
$\reel$-bilinear Lie bracket $\{.,.\}$ on $C^\infty(M)$ satisfying
the Leibniz rule
$$\{f,gh\}=\{f,g\}h+g\{f,h\},\quad\mbox{for all}\; f,g,h\in
C^\infty(M).$$For a function $f\in C^\infty(M)$, the
derivation $H_f=\{f,.\}$ is called the {\bf hamiltonian vector
field} of $f$. If $H_f=0$, we call $f$ a {\bf Casimir function}.
It follows from the Leibniz rule that there exists a bivector field
$\pi\in\Ga(\wedge^2TM)$ such that
$$\{f,g\}=\pi(df,dg);$$the Jacobi identity for $\{.,.\}$ is
equivalent to the condition $[\pi,\pi]=0$, where $[.,.]$ is the
Schouten-Nijenhuis bracket, see e.g. $[13]$.

In local coordinates $(x_1,\ldots,x_n)$ the tensor $\pi$ is
determined by the matrix
$$\pi_{ij}(x)=\{x_i,x_j\}.$$
The rank of this matrix is called the rank of $\pi$ at $x$. A
Poisson structure is called {\bf regular} if the rank of $\pi$ is
constant on $M$.
 If
this matrix is invertible a each $x$, then $\pi$ is called
nondegenerate or {\bf symplectic}. In this case, the local
matrices $(\om_{ij})=(-\pi_{ij})^{-1}$ define a global 2-form
$\om\in\Om^2(M)$, and the condition $[\pi,\pi]=0$ is equivalent to
$d\om=0$.

We denote by  $\an:T^*M\too TM$
 the anchor map given by $\be(\an(\al))=\pi(\al,\be),$
and by $[\;,\;]_\pi$  the Koszul  bracket given by
$$[\al,\be]_\pi=
 L_{\pi_{\#}(\al)}\be-L_{\pi_{\#}(\be)}\al-d(\pi(\al,\be)),\quad \al,
 \be\in \Om^1(M).$$

The distribution $Im\an$ is integrable and defines a singular
foliation $\F$. The leaves of $\F$ are symplectic immersed
submanifolds of $M$. The  foliation $\F$ is called the {\bf
symplectic foliation} associated to the Poisson structure.

The anchor map and the Koszul bracket define the Lie algebroid
structure  associated to $\pi$. The Poisson structure is called
{\bf integrable} if this Lie algebroid structure integrates to a
Lie groupoid structure. (For a detailed explanation of the
integrability of Lie algebroids, see $[5]$, $[6]$ etc.)

\bigskip

The main tool to be used extensively  in this paper
 is the
metric contravariant connection  associated canonically to a
couple of a Riemannian metric and a Poisson tensor. General
contravariant connections associated to a Poisson structure have
recently turned out to be useful in several areas of Poisson
geometry. Contravariant connections were defined by Vaismann
$[13]$ and were analyzed in detail by Fernandes $[7]$. This notion
appears extensively  in the context of noncommutative deformations
(see
 $[12]$, $[8]$ and $[9]$).
One can consult $[7]$ for the general properties
 of contravariant connections.

Let $(M,\pi)$ be a Poisson manifold and $V\stackrel{p}\too M$ a
vector bundle over $M$. A {\bf contravariant connection} on $V$ with
respect to $\pi$ is a map $\D:\Om^1(M)\times
\Ga(M,V)\too\Ga(M,V)$, $(\al,s)\mapsto\D_\al s$ satisfying the
following properties:\begin{enumerate} \item $\D_\al s$ is linear
over $C^\infty(M)$ in $\al$:
$$\D_{f\al_1+h\al_2}s=f\D_{\al_1} s+h\D_{\al_2} s,\quad f,h\in C^\infty(M);$$
\item $\D_\al s$ is linear over $\reel$ in $s$:
$$\D_{\al}(as_1+bs_2)=a\D_\al s_1+b\D_\al s_2,\quad a,b\in \reel;$$
\item $\D$ satisfies the following product rule:
$$\D_\al(fs)=f\D_\al s+\an(\al)(f)s,\quad f\in
C^\infty(M).$$\end{enumerate}

The curvature of a contravariant connection $\D$ is formally
identical to the usual definition
$$K(\al,\be)=\D_\al\D_\be-\D_\be\D_\al-\D_{[\al,\be]_\pi}.$$We
 call $\D$ flat if $K$ vanishes identically.

If $V=T^*M$, one can define the torsion $T$ of $\D$ by
$$T(\al,\be)=\D_\al\be-\D_\be\al-[\al,\be]_\pi.$$

A contravariant connection $\D$  is called  a {\bf
$\F$-connection} if it satisfies the following property
$$\an(\al)=0\qquad\Rightarrow \qquad\D_\al=0;$$
$\D$ will be called  a {\bf $\F^{reg}$-connection} if $\D$ is a
$\F$-connection on the regular open set where the rank of $\pi$ is
locally constant.

 $\F$-connections were introduced by Fernandes in $[7]$; they will appear
 extensively
 in this
 paper.\bigskip

{\bf Remark.}\begin{enumerate} \item The definition of a
contravariant connection is similar to the definition of an
ordinary (covariant) connection, except that cotangent vectors
have taken the place of tangent vectors. So one can translate many
definitions, identities and proof for covariant connections to
contravariant connections simply by exchanging the roles of
tangent and cotangent vectors and replacing Lie Bracket with
Koszul bracket. Nevertheless, Fernandes pointed out in $[7]$ that
some classical results on covariant connections are not true any
more for general contravariant connections and, however, those
results remain  valid for $\F$-connections.

 \end{enumerate}

For a couple of a Poisson tensor and a Riemannian metric there
exists an  unique torsion-free contravariant connection which
preserves the metric. It appeared first in $[1]$ and, recently,
has turned out to be useful in    the context of noncommutative
deformations (see $[8]$ and $[9]$).

 Let $M$ be a pseudo-Rimannian
manifold and $\pi$ a Poisson tensor on $M$. We denote by $<\;,\;>$
the  metric when it measures the length  of 1-forms.

{\bf The metric contravariant connection} associated  to
$(\pi,<\;,\;>)$ is the unique contravariant connection $\D$ with
respect  to $\pi$ such that: \begin{enumerate}\item the  metric
$<,>$ is parallel with respect to $\D$ i.e.,
$$\pi_{\#}(\al).<\be,\ga>=<\D_\al\be,\ga>+<\be,\D_\al\ga>;$$
\item $\D$ is torsion-free.
\end{enumerate}
 One can define $\D$ by the Koszul formula
\begin{eqnarray*}
2<\D_\al\be,\ga>&=&\pi_{\#}(\al).<\be,\ga>+\pi_{\#}(\be).<\al,\ga>-
\pi_{\#}(\ga).<\al,\be>\\
&+&<[\ga,\al]_\pi,\be>+<[\ga,\be]_\pi,\al>+<[\al,\be]_\pi,\ga>.
\end{eqnarray*}

In $[2]$ and $[3]$ we sited  and studied the following notion
of compatibility between a
Riemannian metric and a
 Poisson tensor.
A Poisson tensor $\pi$ on a Riemannian manifold $(M,<,>)$ is
called compatible with the metric
   if:\bigskip

 $(\p1)$ The Poisson tensor $\pi$ is parallel with respect to  the metric
 contravariant
 connection $\D$ associated to $(\pi,<,>)$ i.e., $\D\pi=0$.\bigskip

    We pointed out in $[2]$ and $[3]$ that  a Poisson tensor $\pi$ for  which
    $(\p1)$ holds satisfies the following properties:\bigskip

  $(\p2)$ The divergence of $\pi$ with respect to the Levi-Civita connection
  vanishes.\bigskip

  $(\p3)$ For any open subset $U$ in $M$ and for any function
   $f\in C^\infty(U)$ such that
  $i_{df}\pi=0$, the gradient  vector field $\na f$ preserves $\pi$ i.e.,
  $L_{\na f}\pi=0$.
  \bigskip

  In Section 3, we will recall the definition of the divergence of a
  multi-vector field with
  respect to a covariant connection.

  Note that $(\p2)$ is equivalent to the fact that the Riemannian density is
  invariant by
  any hamiltonian vector field and then $\pi$ is an unimodular Poisson
  structure (see $[15]$).
  \bigskip

  One can remark that the properties $(\p2)$ and $(\p3)$ are not specific
  to bivector fields
   and make sense for  any  multi-vector field  on a Riemannian manifold.
   In particular,
   what can one  tell about a vector field on a Riemannian manifold which
   satisfies $(\p2)$ and
   $(\p3)$? The following proposition gives an answer to this question.
   \begin{pr}
A  vector field $X$ on a Riemannian manifold $(M,g)$ is a Killing vector
field if and
only if the following assertions hold:
\begin{enumerate}\item The divergence of $X$ with respect to the Levi-Civita
connection vanishes. \item For any open set $U\subset M$ and any
function $f\in C^\infty(U)$  such that $X(f)=0$, $[X,\na f]=0$
$(\na f$ is the gradient field of $f$ given by $g(\na
f,Y)=df(Y))$.
\end{enumerate}\end{pr}

{\bf Proof.} Suppose that $X$ is a Killing vector field. It is a
classical result on Killing vector fields that the divergence of $X$ vanishes.
On other hand, for any function $f$ such that $X(f)=0$ and for any vector field
$Y$, we have
$$ 0=L_Xg(Y,\na f)=X.Y(f)-[X,Y](f)-g(Y,[X,\na f])=-g(Y,[X,\na f]),$$ and then
$[X,\na f]=0.$

Conversely, suppose that $X$ is a  vector field which satisfies
the two assertions above.  We will show that $L_Xg$ vanishes on
the dense open set $U_1\cup U_2$ where $U_1=\{m\in M|\;
X(m)\not=0\}$ and $U_2$ is the interior of $\{m\in M|\; X(m)=0\}$.

The vector field $X$ vanishes on $U_2$ and hence $L_Xg(m)=0$ for
any $m\in U_2$.

Let $m$ be a point in $U_1$. Choose a local coordinates
$(x_1,\ldots,x_n)$ such that
 $X=\di\frac{\partial }{\partial x_1}$. The functions $x_2,\ldots,x_n$ satisfy
  $X(x_i)=0$
 and then $[X,\na x_i]=0$ for $i=2,\ldots,n$.

 Let us compute $L_Xg$ on the local frame $(X,\na x_2,\ldots,\na x_n)$. We
  have, for $i=2,\ldots, n$,
\begin{eqnarray*}
L_Xg(X,\na x_i)&=&X.X(x_i)=0.\\
L_Xg(\na x_i,\na x_j)&=&X.\na x_i(x_j)=\na
x_i.X(x_j)=0.\end{eqnarray*} To conclude, it remains to show that
$L_Xg(X,X)=X.g(X,X)=0$. Indeed, we put $E_1=\frac{X}{|X|}$ and we
orthonormalize  $(\na x_2,\ldots,\na x_n)$ to get an orthonormal
frame $(E_1,\ldots,E_n)$. Note that, from $X.g(\na x_i,\na
x_j)=0$, we get $[E_i,X]=0$ for $i=2,\ldots,n$.

The vanishing of the divergence of $X$ gives
$$0=\sum_{i=1}^ng(\na_{E_i}X,E_i),$$where $\na$ is the Levi-Civita connection
associated to the metric. Since , for $i=2,\ldots,n$,
$$g(\na_{E_i}X,E_i)=g(\na_X E_i,E_i)=\frac12X.g(E_i,E_i)=0,$$  we get
$$0=-\frac1{2|X|^2}X.g(X,X),$$which completes the proof.\hskip1cm q.e.d.\bigskip

 According to this proposition, it is natural to put the following definition.

 \begin{Def} A Poisson structure  on a Riemannian manifold $(M,g)$
 will be called a {\bf Killing-Poisson structure}
 if the associated bivector field $\pi$ is a Killing-Poisson tensor i.e.,
 \begin{enumerate} \item the divergence of $\pi$ with respect to the
 Levi-Civita connection
 vanishes,\item for any open set $U\subset M$ and any  function
 $f\in C^\infty(U)$ such that $f$ is a Casimir function, $$L_{\na
 f}\pi=0.$$\end{enumerate}
 \end{Def}

 {\bf Remark.}\begin{enumerate}\item Let $\pi$ be an invertible  Poisson
 tensor on a
Riemannian $2n$-manifold $(M,g)$.  The Riemannian volume $\mu_g$
satisfies $\mu_g=f\wedge^n\om$ where $\om$ is the symplectic form associated
 to $\pi$ and
$f$ a  function. It is easy to check that  $\pi$ is a
Killing-Poisson tensor if and only if $df=0$. Hence, for any
invertible Poisson tensor $\pi$ on a manifold $M$ and for any
Riemannian metric $g$ on $M$, there exists a function $\phi$ such
that $\pi$ is a Killing-Poisson tensor with respect to
$e^{\phi}g$. \item One can check easily that a bivector field on a
Riemannian manifold which is parallel with respect to the
Levi-Civita connection is a Killing-Poisson tensor.
\end{enumerate}

 This paper is organized as follows. In Section 2, we show that $(\p3)$ is
  equivalent to:\bigskip

 $(\p3')$ the metric contravariant connection associated to $(\pi,<,>)$ is
 a $\F^{reg}$-connection.\bigskip

 We show  that a Poisson structure  for which there exists a Riemannian
 metric
 such that $(\p3')$ holds possesses
  the following
properties.

\begin{th} Let $(M,g)$ be a Riemannian manifold endowed with a regular
Poisson tensor $\pi$ such that the metric contravariant connection
associated to $(\pi,g)$ is a $\F$-connection. Then:
 \begin{enumerate} \item the metric $g$ is a bundle-like metric for the
 symplectic foliation,
 \item for any open subset $U$ in $M$ and for any  vector field $X$ on $U$
 which is perpendicular to the symplectic foliation,  $X$ is a foliated
 vector field  if and only if $X$ is a Poisson vector field.\end{enumerate}\end{th}

\begin{th} Let $(M,\pi)$ be a  regular Poisson manifold. If
there exists a Riemannian metric such that the corresponding
metric contravariant connection is a $\F$-connection, then the Lie
algebroid structure    associated to $\pi$ is integrable.\end{th}

Note that the conclusions of Theorem 1.1 and Theorem 1.2 hold, in
particular, for any regular Killing-Poisson structure.\bigskip

In Section 3, we give some general properties of Killing-Poisson structures
and show the following results.
\begin{th} Let $M$ be a Riemannian manifold endowed with a Poisson tensor
such that $\D\pi=0$. Then $\pi$ is a Killing-Poisson tensor.\end{th}

This theorem shows that the class of Killing-Poisson structures is
large. For instance in $[4]$, the author showed that the  dual
$\G^*$ of a Lie algebra $\G$ carries a Riemannian metric
compatible with the canonical linear Poisson structure if and only
if the Lie algebra is a semi-direct product of an abelian Lie
algebra and an abelian ideal.

\begin{th} Let $\G$ be a Lie algebra and $r\in\wedge^2\G$ an unimodular solution
of the classical Yang-Baxter equation. Let $\Ga:\G\too\X(M)$ be a
locally free action  of $\G$ on a Riemannian manifold
$(M,g)$ such that, for any $u\in\G$, $\Ga(u)$ is a Killing vector field.
Then $\pi:=\Ga(r)$ is a
 Killing-Poisson tensor.\end{th}
\begin{co} Any unimodular left-invariant Poisson tensor on a Lie group $G$ is a
Killing-Poisson tensor with respect to any right-invariant
Riemannian metric  on $G$.\end{co} A left-invariant Poisson tensor
on a Lie group is unimodular if its value at the unity is an
unimodular solution of the classical Yang-Baxter equation.

Unimodular solutions of the classical Yang-Baxter equation will be
defined  in Section 3.

In Section 4, we study the exactness of Killing-Poisson tensors
and show  the following result.

\begin{th} Let $\pi$ be a non trivial regular Killing-Poisson tensor on a compact
 Riemannian
 manifold. Then $\pi$ cannot be exact.\end{th}

Section 5 is devoted to the characterization of Killing-Poisson
tensors on  Riemannian 3-manifolds. In particular, we will
characterize all Killing-Poisson tensors on $\reel^3$
 endowed with the Euclidian metric.

 \section{On metric contravariant $\F$-connections}
 In this section, we give some  properties of Poisson   structures
 for which
 there exists a Riemannian metric such that the corresponding metric
 contravariant
 connection is a $\F$-connection, we prove Theorem 1.2 and Theorem 1.3 and
 we show the
 equivalence between $(\p3)$ and $(\p3')$.\bigskip

 Let $M$ be a Riemannian manifold and $\pi$ a Poisson tensor on $M$.
  We denote by $g$ (resp. by $<\;,\;>$) the Riemannian metric when it
  measures the length of vectors  (resp. the length of covectors).
  $\D$ denotes the metric contravariant connection associated to
  $(\pi,g)$, and
  $\#:T^*M\too TM$ the isomorphism associated to the Riemannian
  metric. We denote by $\F$ the symplectic foliation and, for any open set
  $U\subset M$, $\X(\F,U)^\perp$ denotes the space of orthogonal foliated
  vector fields. A vector field $X\in\X(U)$ belongs to $\X(\F,U)^\perp$ if,
  for any vector field $Y\in\X(U)$ tangent to the foliation $\F$, $[X,Y]$ is
  tangent to $\F$ and $g(X,Y)=0$.

 For any open subset $U\subset M$, we consider
 $Z^0(U)=\{f\in C^\infty(U); i_{df}\pi=0\}$ the space of
  Casimir functions on $U$ and $Z^1(U)=\{\al\in\Om^1(U);
  [\al,\be]_\pi=0\;\forall \be\in \Om^1(U)\}$ the center of the
  Lie algebra $\Om^1(U)$ endowed with the Koszul bracket. One can see
  easily that $Z^1(U)$ is the space of basic 1-forms relative to the
  restriction of $\F$ to $U$ i.e.,
 $$Z^1(U)=\{\al\in\Om^1(U); \an(\al)=0\;\mbox{and}\; i_{\an(\be)}d\al=0\;
 \forall \be\in \Om^1(U)\}.$$

 For $\al\in\Om^1(U)$ and $f,h\in C^\infty(U)$, we have
 \begin{eqnarray*}
 L_{\#(\al)}\pi(df,dh)&=&\#(\al).\pi(df,dh)-\pi(L_{\#(\al)}df,dh)-
 \pi(df,L_{\#(\al)}dh)\\
 &=&\#(\al).\pi(df,dh)+\an(dh).<\al,df>-\an(df).<\al,dh>\\
 &=&\#(\al).\pi(df,dh)+<\D_{dh}\al,df>-<\D_{df}\al,dh>+<\al,\D_{dh}df-
 \D_{df}dh>\\
 &=&<\D_{dh}\al,df>-<\D_{df}\al,dh>.\end{eqnarray*}So we get the formula
 $$L_{\#(\al)}\pi(\be,\ga)=
 <\D_{\ga}\al,\be>-<\D_{\be}\al,\ga>,\quad\al,\be,\ga\in\Om^1(U).\eqno(1)$$

 \begin{pr} Let $M$ be a Riemannian manifold, $\pi$ a Poisson tensor on $M$
 and $U$ an open set such that $\D$ is a $\F$-connection on $U$. Then, we have:
 \begin{enumerate}\item for any $\al,\be\in\Om^1(U)$,
 $$\al\in Ker\an\quad \Rightarrow\quad D_\be\al\in Ker\an;$$
 \item for any $\al,\be\in\Om^1(U)$
 $$\al\in Ker\an^\perp\quad \Rightarrow\quad D_\be\al\in Ker\an^\perp,$$
 where $Ker\an^\perp$ is the orthogonal of $Ker\an$.\end{enumerate}\end{pr}

 {\bf Proof. } 1. If $\an(\al)=0$, $\D_\al\be=0$ and then
  $\an(D_\be\al)=\an([\be,\al]_\pi)=[\an(\al),\an(\be)]=0$.

 2. Suppose that $\al\in Ker\an^\perp$ and $\an(\ga)=0$. We have
 $$<D_\be\al,\ga>=\an(\be).<\al,\ga>-<\D_\be\ga,\al>=0$$since
 $\an(\D_\be\ga)=0$ according to 1.\hskip1cm q.e.d.

 \begin{pr} Let $M$ be a Riemannian manifold, $\pi$ a Poisson tensor on $M$
 and $U$ an open set such that $\D$ is a $\F$-connection on $U$. Then
 \begin{eqnarray*}
 Z^1(U)&=&\{\al\in\Om^1(U);\an(\al)=0\;\mbox{and}\; \D\al=0\},\\
 &=&\{\al\in\Om^1(U);\an(\al)=0\;\mbox{and}\; L_{\#(\al)}\pi=0\}.
 \end{eqnarray*}\end{pr}
 {\bf Proof.} The first equality comes from the fact that, if $\D$ is
 a $\F$-connection and $\al\in\Om^1(M)$ such that $\an(\al)=0$, we
 have $[\al,\be]_\pi=-D_\be\al$ for any $\be\in \Om^1(U)$.

 From $(1)$, we have $\D\al=0$ implies $L_{\#(\al)}\pi=0$ and then
 $Z^1(U)\subset
 \{\al\in\Om^1(U);\an(\al)=0\;\mbox{and}\; L_{\#(\al)}\pi=0\}.$

 Suppose now that $\an(\al)=0$ and $L_{\#(\al)}\pi=0$. From Proposition 2.1,
 we get $\an(\D_\be\al)=0$ for any $\be\in\Om^1(U)$. On other hand, let
 $p\in U$ and $\ga\in T^*_pM$ such that $\an(\ga)=0$, by using $(1)$, we
 get for any $\be\in\Om^1(U)$
 $$<\D_\be\al,\ga>=<\D_\ga\al,\be>=0$$and hence $\D_\be\al\in Ker\an\cap
 Kerp\an^\perp$ which implies  $\D\al=0$.\hskip1cm q.e.d.

 \begin{pr} Let $M$ be a Riemannian manifold, $\pi$ a Poisson tensor on $M$
 and $U$ an open set such that $\D$ is a $\F$-connection on $U$ and the rank
 of $\pi$ is constant on $U$. Then
 $$\X(\F,U)^\perp=\#(Z^1(U)).$$\end{pr}

 {\bf Proof.} If $\al\in Z^1(U)$, $\#(\al)$ is perpendicular to the
 symplectic foliation and, according to Proposition 2.2, $\#(\al)$ is a
 Poisson vector field and then a foliated vector field.

 Let $X\in \X(\F,U)^\perp$. Since the Poisson tensor is regular in $U$,
 there exists an unique $\al\in\Om^1(U)$ such that $\an(\al)=0$ and $X=\#(\al)$.
 The vector field $\#(\al)$ is foliated if and only if, for any
 $f\in C^\infty(U)$ and for
 any $h\in Z^0(U)$, $[\#(\al),H_f](h)=0$. This is equivalent to
 $H_f.<\al,dh>=0$ which is equivalent to $<\D_{df}\al,dh>=0$ since $\D dh=0$.
From Proposition 2.1, we deduce that $\D\al=0$ and hence $\al\in
Z^1(U)$, by Proposition 2.2, which completes the proof.\hskip1cm
q.e.d.
 \bigskip

 Let us prove now Theorem 1.1 and Theorem 1.2.\bigskip

 {\bf Proof of Theorem 1.1.} 1. Recall that  the metric $g$ is  bundle-like
 for the foliation $\F$ if it has the following property: for any open set $U$
 of $M$ and for all vector fields $X,Y\in\X(\F,U)^\perp$, the function $g(X,Y)$
 is a Casimir function.

 Let $X,Y\in\X(\F,U)^\perp$. According to Proposition 2.3, $X=\#(\al)$ and
 $Y=\#(\be)$ where $\al,\be\in Z^1(U)$. According to Proposition 2.2,
 $\D\al=\D\be=0$ and then, for any $f\in C^\infty(U)$,
 $$H_f.<\al,\be>=<\D_{df}\al,\be>+<\al,\D_{df}\be>=0$$and the assertion follows.

 2. The result is a consequence of Proposition 2.2 and Proposition 2.3.
 \hskip1cm q.e.d.\bigskip

 {\bf Proof of Theorem 1.2.} Let $M$ be a Riemannian manifold and $\pi$ a
 regular Poisson
tensor on $M$ such that the metric contravariant connection $\D$
is a $\F$-connection. For any symplectic leaf $S$, we have
$$T^*M_{|S}=Ker\pi_{|S}\oplus Ker\pi^{\perp}_{|S},$$ where $Ker\pi$ denotes the
 kernel of $\an$ and $Ker\pi^\perp$ its orthogonal. For any local
 sections $\al$ and $\be$ of $T^*M_{|S}$, we define $[\al,\be]_\pi$ by
$$[\al,\be]_\pi={[\wi\al,\wi\be]_\pi}_{|S},$$where $\wi\al$ and $\wi\be$ are
two extensions of $\al$ and $\be$. In a similar way, we define
$\D_\al\be$.

If $\al,\be$ are two local sections of $Ker\pi^{\perp}_{|S}$, from
Proposition 2.1, $\D_\al\be$ and $\D_\be\al$ are  sections of
$Ker\pi^{\perp}_{|S}$ and hence $[\al,\be]_\pi$ is also a section
of $Ker\pi^{\perp}_{|S}$. We deduce that the inverse
$\an^{-1}:TS\too Ker\pi^{\perp}_{|S}$ of the isomorphism
$\an:Ker\pi^{\perp}_{|S}\too TS$   satisfies
$$[\an^{-1}(X),\an^{-1}(Y)]_{\pi}=\an^{-1}\left([X,Y]\right)\quad\forall
X,Y\in\Ga(TS).$$ Hence $\an^{-1}:TS\too Ker\pi^{\perp}_{|S}$ is a
splitting of the anchor map which is compatible with the Lie
bracket. According to $[5]$ Corollary 5.2 (iii), we get the
 result.\hskip1cm q.e.d.
\bigskip

Let us give now some characterizations of metric contravariant
$\F^{reg}$-connections.

\begin{pr} Let $M$ be a Riemannian manifold and $\pi$ a Poisson tensor on $M$.
 Then the following assertions are equivalent:
\begin{enumerate}\item For any open set $U\subset M$ and for any  function
$f\in Z^0(U)$, $\D df=0$. \item For any open set $U\subset M$ and
for any $\al\in Z^1(U)$, $\D\al=0.$ \item For any open set
$U\subset M$ and for any  function $f\in Z^0(U)$, the gradient
vector field $\na f$ is a Poisson vector field. \item For any open
set $U\subset M$ and for any $\al\in Z^1(U)$, $\#(\al)$ is a
Poisson vector field. \item The metric contravariant connection
$\D$ associated to the metric and $\pi$ is a
$\F^{reg}$-connection.\end{enumerate}\end{pr}

{\bf Proof.} We denote by ${\cal O}^{reg}$ the dense open set
where the rank of $\pi$ is locally constant. Suppose that $\D$ is
a $\F^{reg}$-connection. Let $U\subset M$ be an open set and $f\in
Z^{0}(U)$. We have $df\in Z^1(U)$ and then, according to
Proposition 2.2, $\D df$ vanishes in $U\cap{\cal O}^{reg}$ and
hence in $U$. Moreover, from (1), $L_{\na f}\pi=0$. We have shown
that $5\Rightarrow 1\Rightarrow 3$. In a  similar way, we have
$5\Rightarrow 2\Rightarrow 4.$ On other hand, we have obviously
$4\Rightarrow 3.$ We will establish that
 $3\Rightarrow 5$ and the proposition will follow.

Suppose that $3$ holds. Let $p$ be a regular point of $\pi$ and $\al\in T^*_pM$
 such that $\an(\al)=0.$ There exists an open set $U$ and $f\in Z^0(U)$ such
 that $df(p)=\al$. For $\be,\ga\in\Om^1(U)$, we have
\begin{eqnarray*}
<\D_\al\be,\ga>&=&<\D_{df}\be,\ga>(p)=<\D_\be df,\ga>(p)=<\D_\ga df,\be>(p)\\
&=&<\D_{df}\ga,\be>(p)=-<\D_{df}\be,\ga>(p)=-<\D_\al\be,\ga>\end{eqnarray*}and
hence $\D_\al\be=0$ and the implication follows. \hskip1cm q.e.d.
\section{Killing-Poisson structures}

In this section, we recall the definition of the divergence of a
multi-vector field with respect to a covariant connection and we
prove Theorem 1.3 and Theorem 1.4.\bigskip

 Let $M$ be a smooth manifold and $Q$ a $q$-vector field, that
is, $Q\in\Ga(\wedge^q TM)$ $(q\geq0)$. Let
$$c: \Om^1(M)\times\Ga(\wedge^q TM)\too\Ga(\wedge^{q-1} TM)$$denote the
contraction.

Let $\na$ be a covariant connection on $M$. Then the $(q-1)$-vector field
$div_\na Q$ given by
$$div_\na(Q)=c(\na P)$$ is called the divergence of $Q$ associated with $\na$.

It is shown that if $\na$ is the Levi-Civita connection on an
orientable Riemannian manifold and $Q=X$ is a vector field,
$div_\na X$ is the usual divergence of $X$ with respect to the
Riemannian volume $\mu_g$, i.e., $L_X\mu_g=(div_\na X)\mu_g$.

Although the divergence of a $q$-vector field depends on the choice of the
connection $\na$, we often omit $\na$ and write $div Q$ for $div_\na Q.$

if $\na$ preserves a volume form $\mu$, we have
$$d(i_Q\mu)=-(-1)^qi_{div Q}\mu.\eqno(2)$$
On other hand, if $\pi$ is a Poisson tensor on $M$, we have for any
$f\in C^\infty(M)$,
$$div\pi(f)=div H_f\eqno(3)$$ and then
$$div\pi=0\Leftrightarrow d(i_{\pi}\mu)=0\Leftrightarrow L_{H_f}\mu=0\;
\forall f\in C^\infty(M).\eqno(4)$$

So we get the following characterization of Killing-Poisson structures.
\begin{pr} Let $M$ be a Riemannian manifold and $\pi$ a Poisson tensor on $M$.
Then $\pi$ is a Killing-Poisson tensor if and only if the
following assertions hold:
\begin{enumerate}
\item for any $f\in C^\infty(M)$, $L_{H_f}\mu_g=0$ ($\mu_g$ stands of the
Riemannian density),
\item the metric contravariant connection associated to $\pi$ and the
 Riemannian metric is a $\F^{reg}$-connection.\end{enumerate}\end{pr}

Let us prove now Theorem 1.3.\bigskip

{\bf Proof of Theorem 1.3.} In $[2]$ Section 4, we have shown that
$\D \pi=0$ implies that $div\pi=0$.

Let $U$ be an open subset of $M$, $f\in Z^0(U)$ and
$\al\in\Om^1(U)$. The condition $\D\pi=0$ implies that the
conclusions of Proposition 2.1 hold and then     $\an(\D_\al
df)=0$.

Let $p$ be  regular point  in $U$, $V$  a neighborhood  of $p$ where the rank
of $\pi$ is constant and $\ga\in\Om^1(V)$ such that $\an(\ga)=0$.  There exists
$h\in Z^0(V)$ such that $\al(p)=dh(p)+\al_0$ with $\al_0\in Ker\an^\perp$. Hence
\begin{eqnarray*}
<\D_\al df,\ga>(p)&=&<\D_{df}\al,\ga>(p)=-<\al,\D_{df}\ga>(p)\\
&=&-<dh,\D_{df}\ga>(p)+<\al_0,\D_{df}\ga>(p)=0,\end{eqnarray*}where
$<dh,\D_{df}\ga>=0$ from the definition of $\D$ and
$<\al_0,\D_{df}\ga>=0$ from Proposition 2.1. Hence $\D df=0$ and
then $\D$ is an $\F^{reg}$-connection according to Proposition
2.4.\hskip1cm q.e.d.\bigskip

 To prove Theorem 1.4, we recall some facts about solutions of the classical
Yang-Baxter equation.

 Let $(\G,[\;,\;])$ be a Lie algebra. Recall that    a solution of the
 classical Yang-Baxter equation is a bivector $r\in\wedge^2\G$ such that
 $$[r,r]=0,$$where $[r,r]\in\G\wedge\G\wedge\G$ is defined
by
$$[r,r](\al,\be,\ga)=\al([r(\be),r(\ga)])+\be([r(\ga),r(\al)])+
\ga([r(\al),r(\be)]),$$where $r:\G^*\too\G$ denotes also the linear map given
by $\al(r(\be))=r(\al,\be).$

A solution $r$ of the classical Yang-Baxter equation is called unimodular if
the subalgebra $Imr$ is unimodular i.e., for any $u\in Imr$, the trace of the
endomorphism  $ad_u:Imr\too Imr$, $v\mapsto[u,v]$, vanishes.

Let $\G\stackrel{\Ga}\too\X(M)$ be an action of a Lie algebra $\G$ on a
manifold $M$
 i.e., a morphism of Lie algebras
  from $\G$ to the Lie algebra of vector fields
 on $M$. Any solution $r$ of the classical Yang-Baxter equation defines a
 Poisson tensor $\Ga(r)$ on $M$.

 Let us prove now Theorem 1.4.
 \bigskip

 {\bf Proof of Theorem 1.4.}

 {\bf 1. The divergence of $\Ga(r)$ vanishes}\bigskip

 The solution $r$ defines on $Imr$ a 2-form $\om_r$ by
$$\om_r(u,v)=r(r^{-1}(u),r^{-1}(v)),$$where $r^{-1}(u)$ denotes  any antecedent
of $u$. The 2-form $\om_r$ is nondegenerate and symplectic i.e.,
$$\om_r(u,[v,w])+\om_r(v,[w,u])+\om_r(w,[u,v])=0,\quad u,v,w\in Imr.$$ There
exists a basis
$
(e_1,\ldots,e_n,f_1,\ldots,f_n)$ of $Imr$ such that the symplectic
form $\om_r$ is given by
$$\om_r=\sum_{i=1}^ne_i^*\wedge f_i^*.$$ Since $Imr$ is
unimodular, then for any $z\in Imr$, the trace of $ad_z$ is zero
which  is equivalent to
$$\sum_{i=1}^n\left(\om_r([z,e_i],f_i)+\om_r(e_i,[z,f_i])\right)=0.$$
This relation is equivalent to
$$\sum_{i=1}^n\om_r(z,[e_i,f_i])=0$$and hence to
$$\sum_{i=1}^n[e_i,f_i]=0.$$Denotes by $\mu_g$ the Riemannian volume. Since
$L_{\Ga(e_i)}\mu_g=L_{\Ga(f_i)}\mu_g=0$ for $i=1,\ldots,n$, we have
\begin{eqnarray*}
d(i_{\pi_r}\mu_g)&=&d\left(\sum_{i=1}^ni_{\Ga(e_i)\wedge\Ga(f_i)}\mu_g\right)\\
&=&\sum_{i=1}^n\left(i_{[\Ga(e_i),\Ga(f_i)]}\mu_g-
i_{\Ga(e_i)}L_{\Ga(f_i)}\mu_g-i_{\Ga(f_i)}L_{\Ga(e_i)}\mu_g\right)\\
&=&i_{\Ga(\sum_{i=1}^n[e_i,f_i])}\mu_g=0\end{eqnarray*}and hence $div\Ga(r)=0$.
\bigskip

 {\bf 2. The metric contravariant connection is a $\F^{reg}$-connection}\bigskip

     If
$r=\sum_{i,j}a_{ij}u_i\wedge u_j$, we put, for $\al,\be\in\Om^1(P)$,
$$\D^{r}_\al\be:=\sum_{i,j}a_{ij}\al(U_i)L_{U_j}\be,$$where  $U_i=\Ga(u_i)$.
We get a map $\D^r:\Om^1(P)\times\Om^1(P)\too\Om^1(P)$ which is
the metric   contravariant connection associated to the Poisson
tensor $\Ga(r)$ and the Riemannian metric. On can check easily
that $\D^r$ is a $\F^{reg}$-connection since the action is locally
free. \hskip1cm q.e.d.

\section{Killing-Poisson tensors on a compact Riemannian manifold cannot be
exact}

A Poisson manifold $(M,\pi)$ is called exact if there exists a
vector field $X$ such that $[X,\pi]=\pi$. The vector field $X$ is
called a Liouville vector field. If $\pi$ comes from a symlectic
form $\om$, the condition $[X,\pi]=\pi$  is equivalent to the
exactness of the form $\om$. Although a compact symplectic
manifold cannot be exact, there do exist compact regular Poisson
manifolds admitting Liouville vector fields (see [10], [11] and
[14]). In this section, we show that, like  symplectic structures
on  compact  manifolds, a non trivial regular Killing-Poisson
tensor on a compact Riemannian manifold cannot be exact.

We begin by giving some general results on unimodular exact
Poisson tensors. Recall that a  Poisson manifold $(M,\pi)$ is
called unimodular if there exists a volume form $\mu$ on $M$ such
that  any hamiltonian vector field $H_f$ preserves the volume form
$\mu$ i.e., $L_{H_f}\mu=0$.

Let  $(M,\pi)$ be a Poisson manifold and $X$ a Liouville vector
field on $M$. We have obviously,  for any hamiltonian vector field
$H_f$,
$$[X,H_f]=H_f+H_{X(f)},\quad f\in C^\infty(P).\eqno(5)$$

\begin{pr} Let $(M,\pi)$ be an unimodular exact Poisson manifold. Then,
for any invariant volume form $\mu$, any Liouville vector field $X$ and for
each $n\in\nat$, we have
$$L_X\left(i_{\wedge^n\pi}\mu\right)=di_Xi_{\wedge^n\pi}\mu=(n+div_\mu X)
i_{\wedge^n\pi}\mu.\eqno(6)$$\end{pr}

{\bf Proof.} We have obviously  $[X,\wedge^n\pi]=n\wedge^n\pi$ and
then
$$
ni_{\wedge^n\pi}\mu=i_{[X,\wedge^n\pi]}\mu =L_X\circ
i_{\wedge^n\pi}\mu-i_{\wedge^n\pi}\circ L_X\mu =L_X\circ
i_{\wedge^n\pi}\mu-div_\mu X(i_{\wedge^n\pi}\mu).$$This gives the
relation.\hskip1cm q.e.d.
\begin{pr} Let $(M,\pi)$ be  a compact
unimodular regular non trivial Poisson manifold. Then, for any
Liouville vector field $X$, there exists a point in $M$ where $X$
is transverse to the symplectic foliation.\end{pr}

{\bf Proof.} Let $(M,\pi)$ be  a compact unimodular regular
Poisson manifold of rank $2q$ $(q>0)$. Suppose that there exists a
Liouville vector field $X$ which is everywhere tangent  to the
symplectic foliation.  This  implies that the multi-vector field
$X\wedge(\wedge^q\pi)$ vanishes identically. From (6), we get
that, for any invariant volume form $\mu$, the form $(q+div_\mu
X)i_{\wedge^q\pi}\mu$ vanishes also. The form $i_{\wedge^q\pi}\mu$
being a transverse volume form to the symplectic foliation, we get
$div_\mu X=-q$ which is a contradiction with $\int_Mdiv_\mu
X\mu=0$.\hskip1cm q.e.d.\bigskip

We are  able now to give a proof of Theorem 1.5.\bigskip

{\bf Proof of Theorem 1.5.} Let $X$ be a Liouville vector field.
The vector field $X$ splits $X=X^t+X^{\perp}$ where $X^t$ is
tangent to the symplectic foliation and $X^\perp$ is perpendicular
to the symplectic foliation. From (5), $X$ is a foliated vector
field which implies that $X^\perp$ is a foliated vector field and
hence a Poisson vector field according to Theorem 2.1. This
implies that $X^t$ is a Liouville vector field which it
contradicts Proposition 4.2.\hskip1cm q.e.d.

\section{Killing-Poisson structures on 3-dimensional Riemannian manifolds}

We begin by giving a description of Killing-Poisson tensors on a
Riemannian surface. Let $S$ be a connected  orientable Riemannian
surface and $\pi$ a Poisson tensor on $S$. If $\mu$ is the
Riemannian volume, we have that $div\pi=0$ if and only if the
function $i_\pi\mu$ is constant. This mean that $\pi$ is either
trivial or symplectic. Hence a non trivial Killing-Poisson tensor
on $S$ is symplectic.

\begin{th} Let $M$ be a 3-dimensional oriented Riemannian manifold and $\pi$ a
bivector field on $M$.   We denote by $\mu$ the  Riemannian volume. Then the
 following assertions are equivalent:
\begin{enumerate}
\item $\pi$ is a Killing-Poisson tensor.
\item $\D\pi=0$.
\item The 1-form $\al=i_\pi\mu$ satisfies:
$$d\al=0\quad\mbox{and}\quad d<\al,\al>+\de(\al)\al=0,$$
where $\de(\al)=-div(\#(\al))$.\end{enumerate}\end{th}

{\bf Proof.} Recall (see [13]) that for any muliti-vector fields
$Q$ and $R$  on $M$, we have
$$i_{[Q,R]}\mu=(-1)^{(|Q|-1)(|R|-1)}\left(i_Qdi_R\mu-(-1)^{|Q|}di_Qi_R\mu
+(-1)^{(|Q|-1)|R|+|Q|}i_Rdi_Q\mu\right)$$where $|Q|$ is the degree
of the  multi-vector field $Q$.

From this relation, we get for $X$  a vector field and $\pi$ is a
bivector field:
$$i_{[\pi,\pi]}\mu=di_{\pi\wedge\pi}\mu-2i_\pi di_\pi\mu,\eqno(8)$$
$$i_{[X,\pi]}\mu=i_Xdi_\pi\mu+di_Xi_\pi\mu-(div_\mu X)i_\pi\mu.\eqno(9)$$

We will prove now the equivalence $1\Leftrightarrow 2$.

We have seen, in Theorem 1.3, that $2\Rightarrow1$.

Conversely, suppose that $\pi$ is a Killing-Poisson tensor. The
1-form $\al=i_\pi\mu$ satisfies  $\an(\al)=0$ and $d\al=0$ and
then $\al\in Z^1(U)$. From Proposition 2.2, we get $\D\al=0$.
Since $\D\mu=0$, we get $\D\pi=0$ which completes the proof of the
equivalence.

We will prove now the equivalence $1\Leftrightarrow 3$.

Suppose that $\pi$ is a Killing-Poisson tensor. We have $d\al=0$
and $\al\in Z^1(U)$ and, from Proposition 2.2, $[\#(\al),\pi]=0$.
Then, from (9), we get the relation
$$ d<\al,\al>+\de(\al)\al=0.$$
Suppose now that 3 holds. From (8), we get   $[\pi,\pi]=0$ and
then $\pi$ is a Poisson tensor, $div\pi=0$ and also
$[\#(\al),\pi]=0$. Let $U$ be an open set and $f\in Z^0(U)$. We
will show that $\na f$ is a Poisson vector field and the result
follows.

We denote by $U^{reg}$  the open set intersection of $U$ with the
dense open set of  regular point of $\pi$. Let $p\in U^{reg}$. We
have two case:

First case: $\al(p)=0$ and then $\pi$ vanishes on a neighborhood
of $p$ and  hence $[\na f,\pi](p)=0$.

Second case: $\al(p)\not=0$ and then $\al$ does not vanish on a
neighborhood of $p$ and hence $df=h\al$ where $h$ is a local
function. Since $\al\in Z^1(M)$, $h$ is a Casimir function and
then $[\na f,\pi]=h[\#(\al),\pi]=0$. This completes the
proof.\hskip1cm q.e.d.

\begin{co} Let $M$ be an oriented Riemannian 3-manifold
 such that $H^1_{dR}(M)=0$ and let $\pi$ be a bivector
field on $M$. The bivector field $\pi$ is a Killing-Poisson tensor
if and  only if there exists $f\in C^\infty(M)$ such that
$i_{\pi}\mu=df$ and
$$d<df,df>+\De(f)df=0$$where $\De$ is the Beltrami-Hodge Laplacian acting
on functions.\end{co}

{\bf Examples.}  According to Corollary 5.1, a bivector field
$\pi$ on $\reel^3$ is a Killing-Poisson tensor with respect to the
Euclidian  metric if and only if
$$\pi=\frac{\partial f}{\partial z}\frac{\partial}{\partial x}\wedge
\frac{\partial}{\partial y}-\frac{\partial f}{\partial y
}\frac{\partial}{\partial x}\wedge \frac{\partial}{\partial z}+
\frac{\partial f}{\partial x}\frac{\partial}{\partial y}\wedge
\frac{\partial}{\partial z}$$where $f\in C^\infty(\reel^3)$
verifies
$$d<df,df>+\De(f)df=0.\eqno(E)$$

The polynomial functions of degree 2 solutions of $(E)$ are
$$
f(x,y,z)=(a+c)x^2+(a+b)y^2+(b+c)z^2
-2\sqrt{bc}xy+2\sqrt{ab}xz+2\sqrt{ac}yz$$ where $a,b,c$ are real
constants with  the same sign. This gives all linear
Killing-Poisson structures on $\reel^3$ endowed with the Euclidian
metric.\bigskip

One can check easily that a function $f(x,y,z)=g(r)$ where
$r=x^2+y^2+z^2$ is a solution of $(E)$ if and only if the function
$g$ satisfies the differential equation
$$2ry''-y'=0.$$Then $g(r)=ar^{\frac32}$ where $a$ is a constant.
Hence, the Poisson tensor
$$\pi=\sqrt{x^2+y^2+z^2}\left(z\frac{\partial}{\partial x}\wedge
\frac{\partial}{\partial y}-y\frac{\partial}{\partial x}\wedge
\frac{\partial}{\partial z}+ x\frac{\partial}{\partial y}\wedge
\frac{\partial}{\partial z}\right)$$is a Killing-Poisson tensor on
$\reel^3$. Remark that $\pi=\sqrt{x^2+y^2+z^2}\pi_{so(3)}$ where
$\pi_{so(3)}$ is the Lie-Poisson structure on the dual of the Lie
algebra $so(3)$.\bigskip

 {\bf Acknowledgements}\bigskip

Thanks to A. El Soufi who invited me at Fran\c cois Rabelais
University at  Tours where this work was written. \bigskip

{\bf References}\bigskip

 [1] {\bf
M. Boucetta, } {\it Compatibilit\'e des structures
pseudo-riemanniennes et des structures de Poisson,}{ C. R. Acad.
Sci. Paris, {\bf t. 333}, S\'erie I, (2001) 763--768.}

[2] {\bf M. Boucetta,  } {\it Poisson manifolds with compatible
pseudo-metric and pseudo-Riemannian Lie algebras,}{ Differential
Geometry and its Applications, {\bf Vol. 20, Issue 3}(2004),
279--291.}

[3] {\bf M. Boucetta, } {\it Riemann-Poisson manifolds and
K\"ahler-Riemann foliations,}{ C. R. Acad. Sci. Paris, {\bf t.
336}, S\'erie I, (2003) 423--428.}

 [4] {\bf M. Boucetta,} {\it  On the Riemann-Lie algebras and
 Riemann-Poisson Lie groups,} J. Lie Theory 15 (2005), no. 1, 183--195.

 [5] {\bf M. Crainic and  R. Fernandes,  }
  {\it Integrability of Lie brackets,} Ann. of Math. (2) 157 (2003),
  no. 2, 575--620.

 [6] {\bf M. Crainic and R. Fernandes, }
  {\it Integrability of Poisson brackets,} J. Differential Geom. 66 (2004),
  no. 1, 71--137.

  [7] {\bf R. Fernandes, } {\it  Connections in Poisson geometry. I.
  Holonomy and invariants,} J. Differential Geom. 54 (2000), no. 2, 303--365.

  [8] {\bf E. Hawkins,} {\it Noncommutative Rigidity}, Commun. Math.
Phys. {\bf 246} (2004) 211-235. math.QA/0211203.

[9] {\bf E. Hawkins ,} {\it The structure of noncommutative
deformations,}

arXiv:math.QA/0504232.

[10] {\bf   G. Hector,  E. Macias  and  M. Saralegi,} {\it Lemme
de Moser feuillet\'e et classification des vari\'et\'es de Poisson
r\'eguli\`eres,} Publ. Mat. {\bf 33} (1989) 423-430.

[11] {\bf    T. Mizutani, } {\it  On exact Poisson manifolds of
dimension 3,}
 Foliations: Geometry and dynamics (Warsaw, 2000), 371-386,
 World Sci. Publishing, River Edge, NJ, 2002.

[12] {\bf N. Reshetikhin, A. Voronov, A. Weinstein,} {\it
Semiquantum Geometry, Algebraic geometry}, J. Math. Sci. {\bf
82(1)} (1996) 3255-3267.

 [13] {\bf I. Vaisman,} {\it Lecture on
the geometry of Poisson manifolds}, Progr. In Math. {\bf Vol.
118}, Birkhausser, Berlin, (1994).

[14] {\bf  A. Weinstein  and P. Xu,} {\it Hochschild cohomology
and characteristic classes for star-products}, Geometry of
differential equations, 177--194, Amer. Math. Soc. Transl. Ser. 2,
186, Amer. Math. Soc., Providence, RI, 1998.

[15] {\bf  A. Weinstein,} {\it  The Modular Automorphism Group of
a Poisson Manifold}, J. Geom. Phys. {\bf23}, (1997) 379-394.

[16] {\bf  Weinstein A.,} {\it  The Local Structure of Poisson Manifolds,} J.
Differential Geometry {\bf18}, (1983) 523-557.

 \end{document}